\documentclass{amsart}
\usepackage{amsmath,amsthm,amsfonts,amsbsy,amssymb,amsmath,graphics,graphicx,mathrsfs}
\usepackage[english]{babel}
\usepackage[matrix,arrow]{xy}
\usepackage{mathptmx}
\usepackage{setspace}
\usepackage{amscd}
\usepackage{epstopdf}
 \theoremstyle{plain}    
 
 \newtheorem{theorem}{Theorem}[section]
 \numberwithin{equation}{section} 
 \numberwithin{figure}{section} 
 \theoremstyle{definition}
 \newtheorem{definition}[theorem]{Definition}
 
 \theoremstyle{plain}    
 
 \newtheorem{proposition}[theorem]{Proposition}
 \theoremstyle{remark}
 \newtheorem{remark}[theorem]{Remark}
 \newtheorem{example}[theorem]{Example}

\theoremstyle{plain}    
 \newtheorem{lemma}[theorem]{Lemma} 
 
 \theoremstyle{definition}
 \theoremstyle{plain}    

\newcommand{\Z}{\mathbb{Z}}

\newcommand{\R}{\mathbb{R}}

\DeclareMathOperator{\Cl}{Cl}

\usepackage{epsfig}

\begin{document}

\title{Markov Theorem for Free Links}

\author{Vassily Olegovich Manturov}
\address{Peoples' Friendship University of Russia, Moscow 119991, Ordjonikidze St., 3}
\email{vomanturov@yandex.ru}
\author{Hang Wang}
\address{Mathematical Sciences Center Room 143, Jin Chun Yuan West Bldg., Tsinghua University, Haidian Dist., Beijing 100084, China.}
\email{hwang@math.tsinghua.edu.cn}

\begin{abstract}
The notion of free link is a generalized notion of virtual link. In the present paper we define the group of free braids, prove the Alexander theorem that all free links can be obtained as closures of free braids and prove a Markov theorem, which gives necessary and sufficient conditions for two free braids to have the same free link closure. Our result is expected to be useful in study the topology invariants for free knots and links.
\end{abstract}

\maketitle




\section{From classical knots to free knots}

Knot theory studies the isotopy classes of smooth embeddings of $S^1$ into three-sphere $S^{3}$ (or, equivalently, to three-space $\R^3$).
A {\em \(classical\) knot} is the image of a smooth embedding of the circle $S^{1}$;
two knots are {\em isotopic}\index{Knots!isotopic}, if one of
them can be transformed into the other by an orientation-preserving diffeomorphism of the
ambient space $S^{3}$. If we
embed a disjoint union of several circles $S^{1}\sqcup \dots \sqcup S^{1}$ in $S^{3}$, then we obtain a {\em classical link};
knots are encoded by their knot diagrams, which are images of smooth immersions of the circle in a plane with an additional structure. 
In the paper, we use generic terms ``links" for both knots and links.

\begin{definition}
A {\em link diagram} is a framed $4$-graph. Each vertex of this graph, also called a {\em crossing
of a link diagram}, is endowed with the structure of an {\em overcrossing} or a {\em
undercrossing}. See Fig. \ref{overunder}.

\begin{figure}
\begin{center}
\begin{picture}(50,50)
\thicklines\put(5,5){\line(1,1){40}} \put(45,5){\line(-1,1){18}}
\put(23,27){\line(-1,1){18}} \put(50,40){overcrossing}
\put(50,10){undercrossing} \thicklines
\end{picture}
\end{center}
\caption{The local structure of a crossing} \label{overunder}
\end{figure}
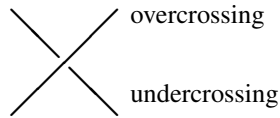

Here, a graph whose vertices have the same degree $k$ is called a {\em
$k$-graph} and a $4$-graph is called {\em framed}\index{Graph!framed} if for every
vertex the four emanating half-edges are split into two pairs of
formally opposite edges. We allow loops and multiple edges.
We restrict ourselves to finite graphs only.
\end{definition}

\begin{figure}
\centering\includegraphics[width=250pt]{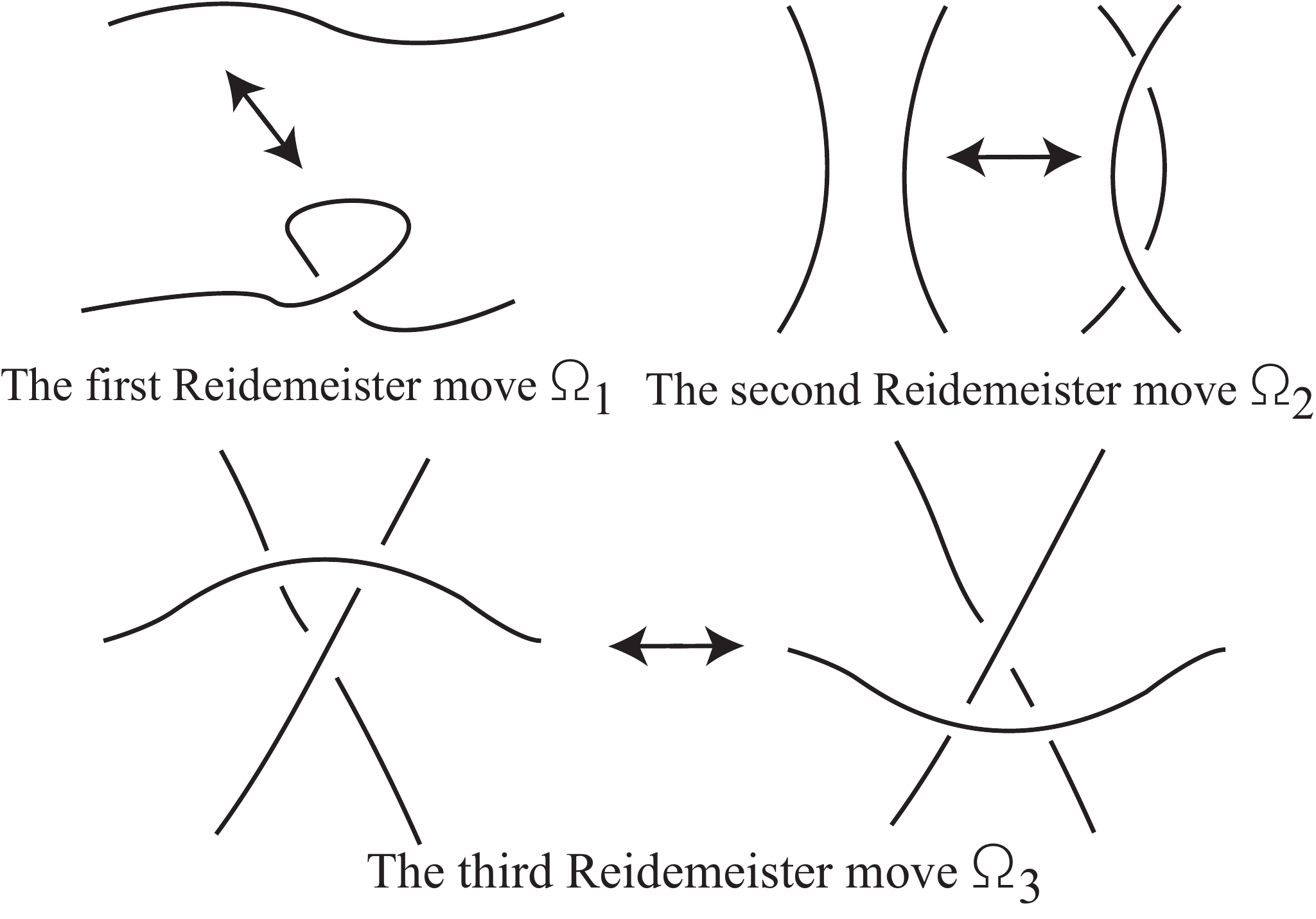}
\caption{Reidemeister moves $\Omega_{1},\Omega_{2},\Omega_{3}$}
\label{rms}
 \end{figure}

It is well known that two link diagrams represent isotopic links if and only if one can be transformed to the other by a sequence of planar isotopies and {\em Reidemeister
moves} \cite{Rei}, see Fig.~\ref{rms}. The Reidemeister theorem allows one to consider isotopy classes of links as
combinatorial objects, which represent equivalence classes of planar
diagrams under Reidemeister moves.

A remarkable generalization of knot theory is virtual knot theory, invented by Louis Kauffman in mid 1990-s \cite{ViK}. Virtual knot theory describes knots and links in thickened surfaces up to isotopy, {\em stabilization} and {\em destabilization}. Here, by stabilization (or destabilization) we mean addition (or removal) of a thickened handle to the thickened surface away from the knot/link. 
A {\em virtual knot} (or, in the case of many components, a {\em virtual link}) represents a natural combinatorial generalization of a classical knot (or link): we introduce a new type of crossing, that is, {\em virtual crossing}, and add new moves, which can be together described as local version of the global {\em detour move} (See Fig. \ref{detour move}), to the list of the Reidemeister moves. 

\begin {definition}
A {\em virtual diagram} is the image of an
immersion of a framed $4$-graph in ${\mathbb R}^2$ with a finite
number of generic projections of edges, that means, the edges are transverse to each other and the intersections point of a crossing has at most $2$-pre-images. 
Here, the immersion respects framing, that is, formally opposite edges remain opposite on the plane.
Each vertex of the graph is endowed with
the classical crossing structure (with a choice for underpass and
overpass specified). The images of vertices with such an additional
structure are called {\em classical crossings}.
Moreover, intersection points of image of edges are called {\em virtual crossings} and
are marked by small circles. 

A {\em virtual link} is an equivalence class of
virtual diagrams modulo planar isotopies and generalized Reidemeister moves. The latter
consist of the usual Reidemeister moves for the classical
crossings and the {\em detour move}, which can be viewed as a replacement of an arc of a virtual link containing only virtual crossings connecting some point $A$ to some other point $B$ of the virtual diagram, by another arc of such sort drawn elsewhere in the plane; all new crossings of the new arc are to be virtual, see Fig \ref{detour move}. 
 \begin{figure}
\centering\includegraphics[width=220pt]{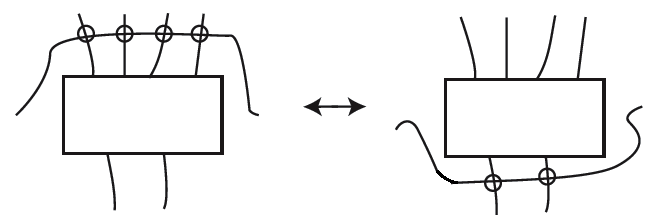}
\caption{A detour move}\label{detour move}
 \end{figure}
 \end {definition}

The detour move can be also expressed as a sequence of ``local detour moves" or ``generalized Reidemeister moves":
 \begin{enumerate}
  \item
{\em Virtual Reidemeister moves}\index{Reidemeister moves!virtual}
$\Omega'_{1}, \Omega'_{2},\Omega'_{3}$, which are obtained from the
classical Reidemeister moves by swapping all classical crossings
participating in the moves for virtual crossings, see Fig.~\ref{vver}.

 \begin{figure}
\centering\includegraphics[width=200pt]{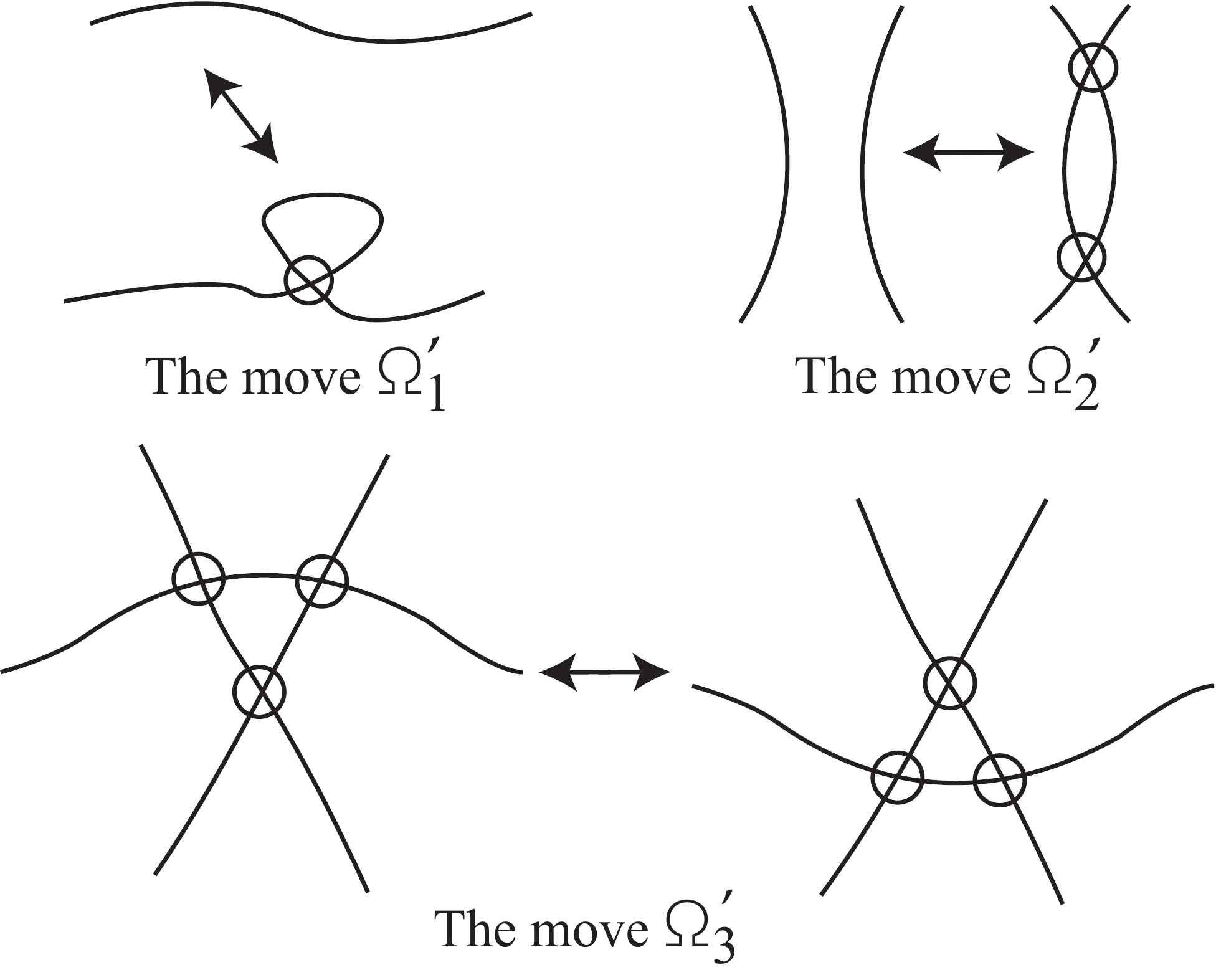}
\caption{Moves $\Omega'_{1},\Omega'_{2},\Omega'_{3}$}\label{vver}
 \end{figure}

 \item
{\em Semivirtual Reidemeister move} $\Omega''_{3}$. Under this move the branch containing
two virtual crossings can slide through a classical crossing,
see Fig.~\ref{semvir}.

 \begin{figure}
 \centering\includegraphics[width=150pt]{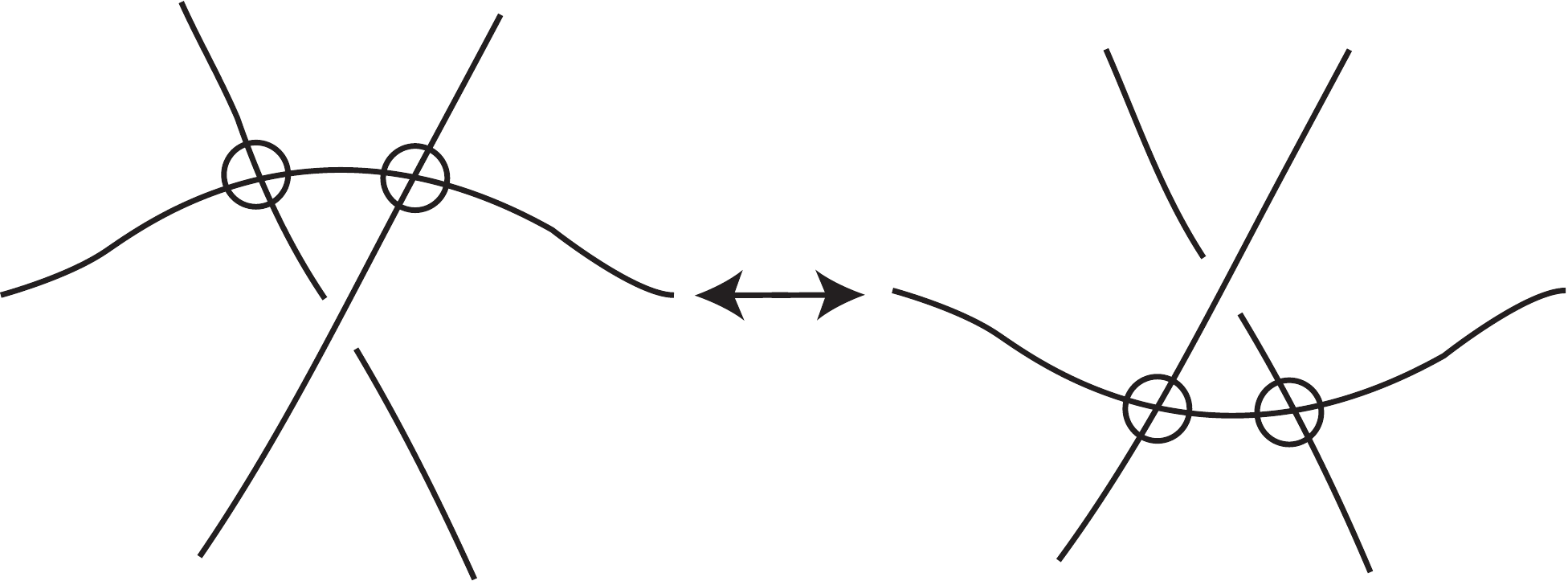}
\caption{The semivirtual move $\Omega''_{3}$}\label{semvir}
\end{figure}
  \end {enumerate}

\begin{figure}
\centering\includegraphics[width=150pt]{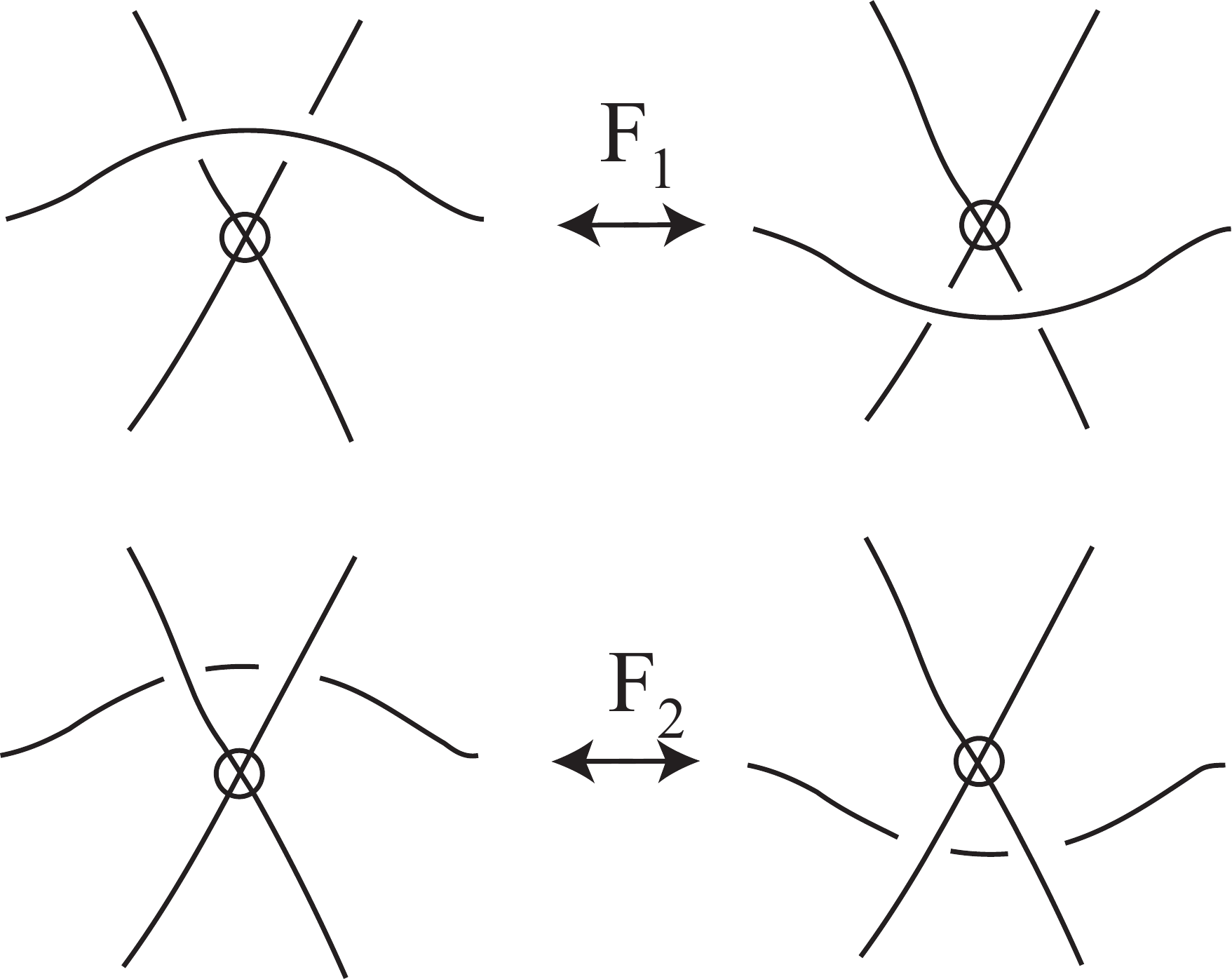}
\caption{Forbidden moves of virtual diagrams}\label{rdmst_forb}
\end{figure}

We note that the \emph{forbidden moves} shown in Fig. \ref{rdmst_forb} are not in the list of generalized Reidemeister moves. Moreover, these moves are not consequences of the generalized Reidemeister moves, see, for example, \cite{KaV, Man1}.

\begin{figure}[t]
\centering\includegraphics[width=350pt]{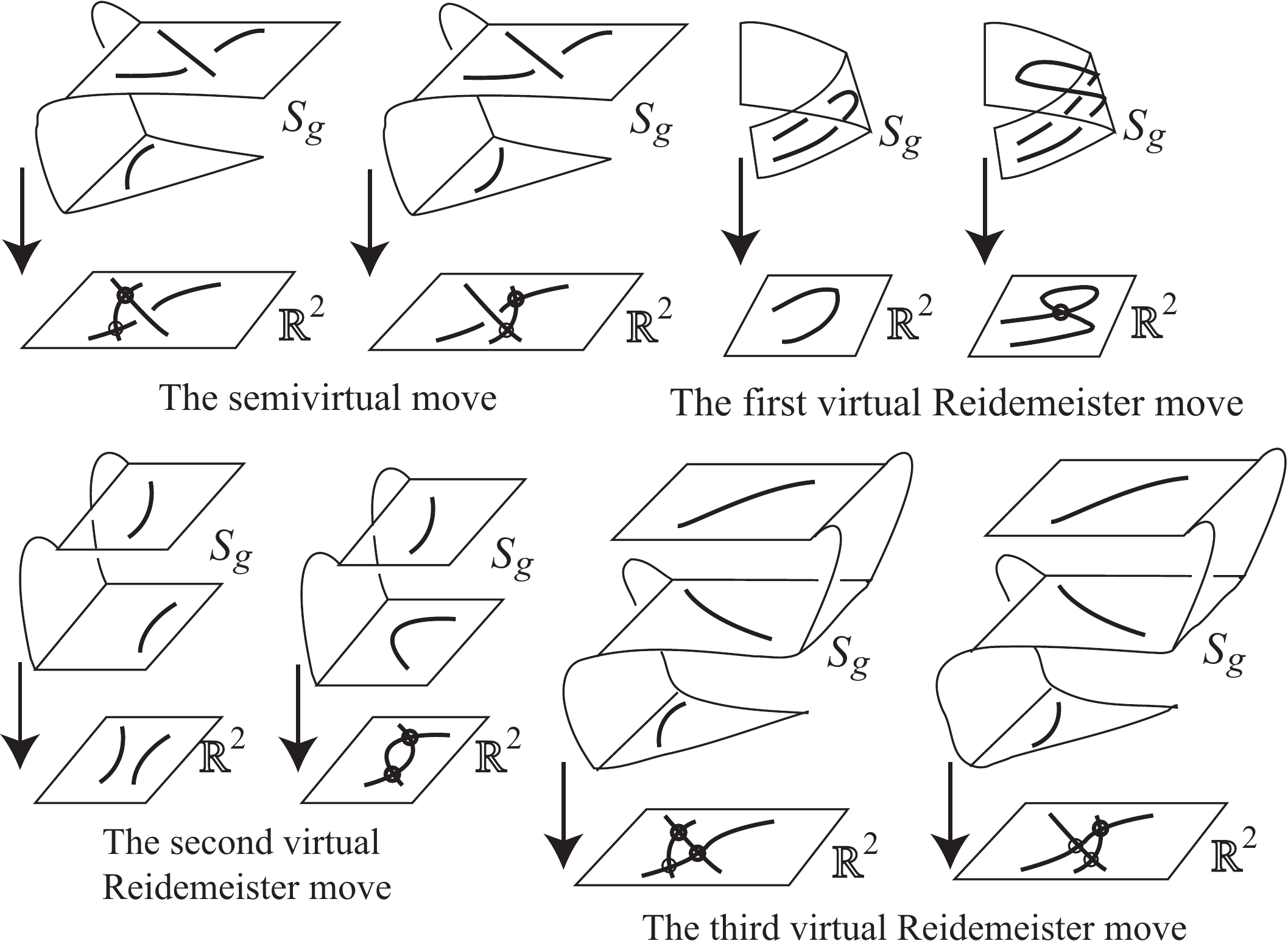}
\caption{Generalized Reidemeister moves and thickened surfaces}
\label{huabao}
 \end{figure}

As mentioned above, virtual links can be realized as links in thickened
oriented surfaces; moreover, thickened surfaces should be considered up to
stabilizations and destabilizations. The Reidemeister moves for diagrams on such a thickened surface $S$ correspond to the classical Reidemeister
moves for virtual diagrams; there are also transformations which do
not change the combinatorial structure of a diagram on $S$, but do
change the combinatorial structure of the projection to $\R^2$: transformations of such sort are described by detour moves. A realization of the detour move by moves on thickened surfaces and
their projections is shown in Fig.~\ref{huabao}.

Classical knot theory admits a fundamental approach to invariants
(for example, skein-invariants, quantum invariants, Vassiliev invariants),
which is based on the following simple observation: every classical knot
can be unknoted by using crossing switches.
For virtual knots, it is not so. Moreover, it is not hard to see that if we
consider equivalence classes of virtual knots by crossings switches, 
we get homotopy classes of curves on surfaces modulo stabilization/destabilization. This operation corresponds to the passage from thickened
surfaces to bare 2-surfaces. This is what is called {\em flat virtual links}, see, for example,~\cite{KK}, on which we could extend invariants of classical knots. More precisely, we give the following definition.

\begin{figure}
\centering\includegraphics[width=270pt]{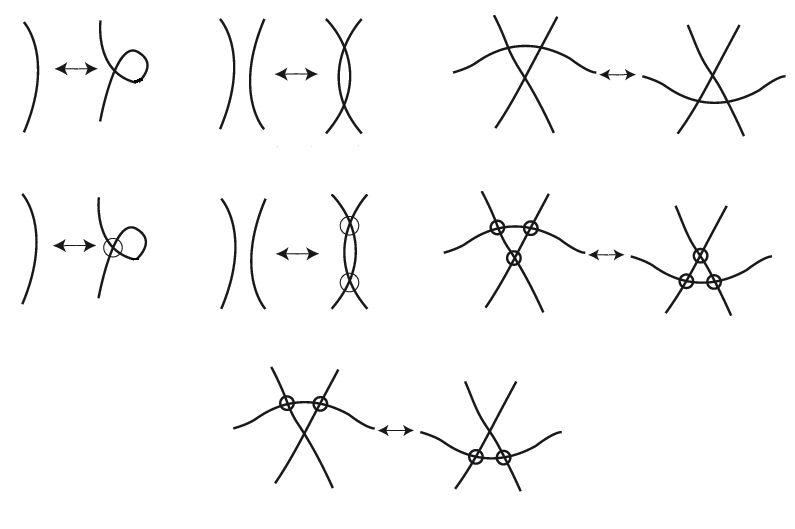}
\caption{Reidemeister moves for flat virtual links} \label{flmoves}
 \end{figure}

\begin{definition}
A \emph{flat virtual diagram} is the image of an immersion of a framed $4$-graph in $\R^2$ with a finite number of generic projections of edges, where the intersection points of images of edges are virtual crossings, represented as crossings with a small circle. The images of the vertices are called \emph{flat crossings}, represented as crossings without decorations.
A {\em flat virtual link} is an equivalence class of flat
virtual diagrams modulo planar isotopies and generalized Reidemeister moves for flat virtual links, which are depicted in Fig.~\ref{flmoves}.
\end{definition}

\begin{figure}
 \centering\includegraphics[width=120pt]{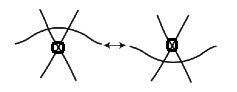}
\caption{Forbidden moves for virtual links and for flat virtual links} \label{fm}
 \end{figure}

It is obvious that the forbidden move for flat virtual links in Fig.~\ref{fm} is not in the list and that it is not a consequence of the moves from the list of Reidemeister moves for flat virtual links.
 
\begin{figure}
\begin{center}
\begin{picture}(100,160)
\thicklines \put(5,95){\line(1,1){50}} \put(55,95){\line(-1,1){50}}
\put(50,120){$\longrightarrow$} \put(50,40){$\longrightarrow$}
\put(5,5){\line(1,1){50}} \put(55,5){\line(-1,1){50}}
\put(30,30){\circle{5}} \thinlines \put(65,100){\line(1,1){20}}
\put(85,120){\line(-1,1){20}} \put(70,95){\line(1,1){20}}
\put(90,115){\line(1,-1){20}} \put(115,100){\line(-1,1){20}}
\put(95,120){\line(1,1){20}} \put(110,145){\line(-1,-1){20}}
\put(90,125){\line(-1,1){20}} \put(65,10){\line(1,1){50}}
\put(70,5){\line(1,1){50}} \put(115,5){\line(-1,1){20}}
\put(120,10){\line(-1,1){20}} \put(65,55){\line(1,-1){20}}
\put(70,60){\line(1,-1){20}}
\end{picture}
\end{center}
\vspace{-0.5cm} \caption{The local structure of $M'$} \label{cr}
\end{figure}
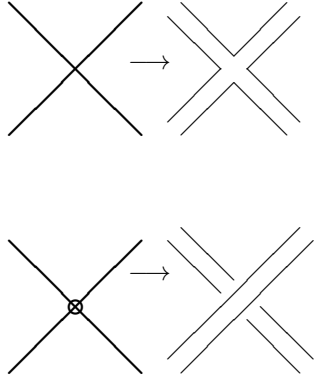

Flat virtual links are realized as homotopy
classes of curves on $2$-surfaces considered up to {\em stabilization} (an addition of a handle to a surface after removing two discs disjoint from our curve) and {\em destabilization} (the inverse operation).  
In fact, we have a natural lifting of flat virtual links to $2$-surfaces. 
It is constructed in two steps \cite{KK}. First of all, having a flat
virtual diagram $L$ we construct a surface with a boundary as
follows. In each flat crossing of a link diagram we place a
cross made of two flat intersecting bands (the upper picture of Fig.~\ref{cr}).
At each virtual crossing we set two nonintersecting bands (the lower picture), cf.~\cite{KK}.
Connecting these crosses and bands by bands (non-overtwisted) along
the arcs of the link we obtain an oriented $2$-dimensional manifold
with boundary. Denote the resulting manifold by $M'$.

One can naturally project the diagram of $L$ to $M'$ in such a way
that arcs of the diagram are projected to middle lines of bands;
herewith flat (classical) correspond to crossings in ``crosses''.
Thus, we obtain a set of closed curves $\delta\subset M'$. Attaching discs
to the boundary components of $M'$, one obtains an orientable surface
$M=M(L)$ without boundary with the set $\delta$ of circles immersed
in it.
This leads us to the following theorem.

 \begin {theorem}[\cite{KK}]
Flat virtual links are equivalence classes of finite sets of curves
in $2$--surfaces up to free homotopy, stabilization and
destabilization.
 \end {theorem}

Let us pass to free knot theory. 
It follows from the definition that flat virtual links are obtained by factorization of virtual links \cite{IMN11}; As for free links, they are proved to be a quotient of flat virtual links; however, it is not known whether they are 
algorithmically recognizable. The reason is that they have a combinatorial definition but there is no obvious geometry behind them, unlike the case of flat virtual links \cite{KK}.
Therefore, free links are interesting and we shall focus on them.

\begin{definition}[\cite{ManOFK}]
A {\em free link} is an equivalence class of framed $4$-valent graphs modulo the following three transformations:\par
\begin{enumerate}
\item {\em The first Reidemeister move} being an addition/removal of a loop, see Fig.~\ref{1r}.
 \begin{figure}
\centering\includegraphics[width=170pt]{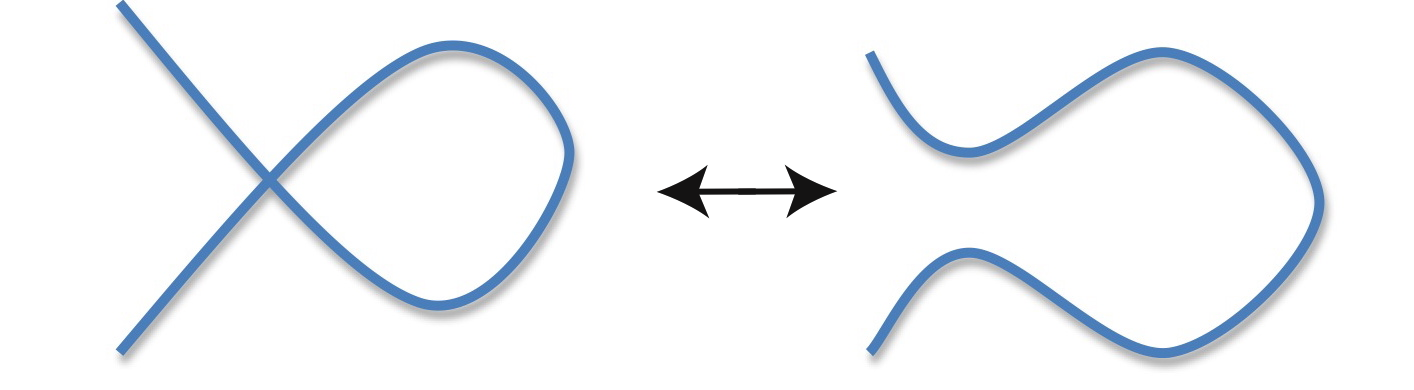} \caption{The first
Reidemeister move for free links} \label{1r}
 \end{figure}
\item {\em The second Reidemeister move} being an addition/removal of a bigon formed by a pair of edges which are adjacent (not opposite) at each of the two vertices, see Fig.~\ref{2r}.
 \begin{figure}
\centering\includegraphics[width=220pt]{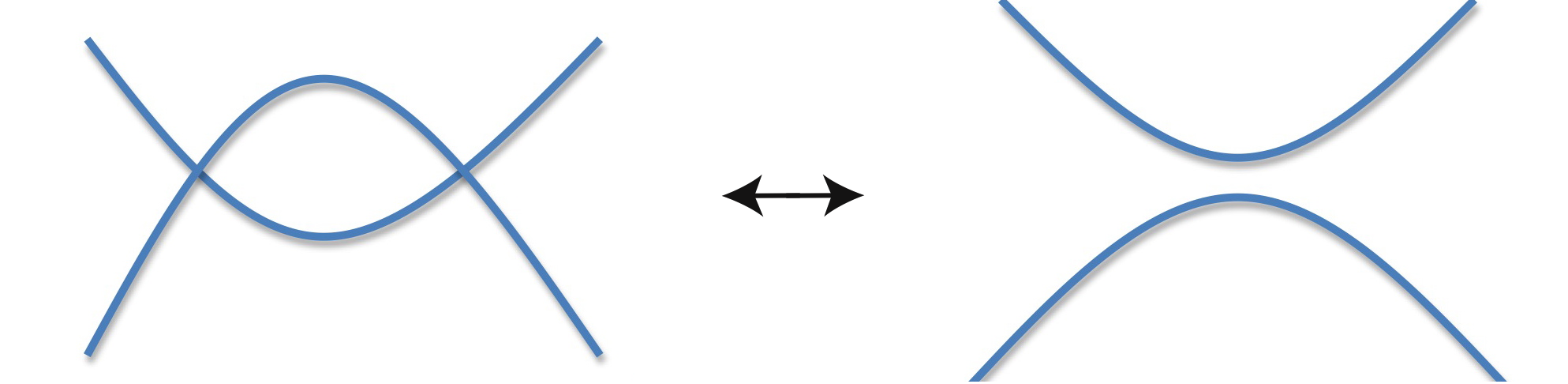} \caption{The second
Reidemeister move for free links} \label{2r}
 \end{figure}
\item {\em The third Reidemeister move} being a triangle move involving three vertices, see
Fig.~\ref{sootvt}.

 \begin{figure}
\centering\includegraphics[width=150pt]{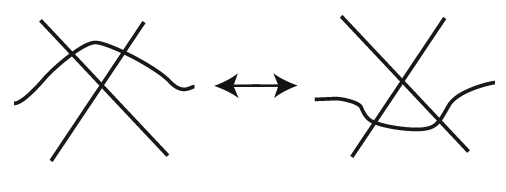} \caption{The
third Reidemeister move for free links} \label{sootvt}
 \end{figure}

\end{enumerate}
Here, every move respects the framing. 
When we project it on a plane, we obtain a {\em free link diagram}, a framed $4$-valent graph where we have flat crossings as well as virtual crossings. 
The framing of a free link diagram is naturally taken from the plane.
\end{definition}

\begin{figure}
\centering\includegraphics[width=150pt]{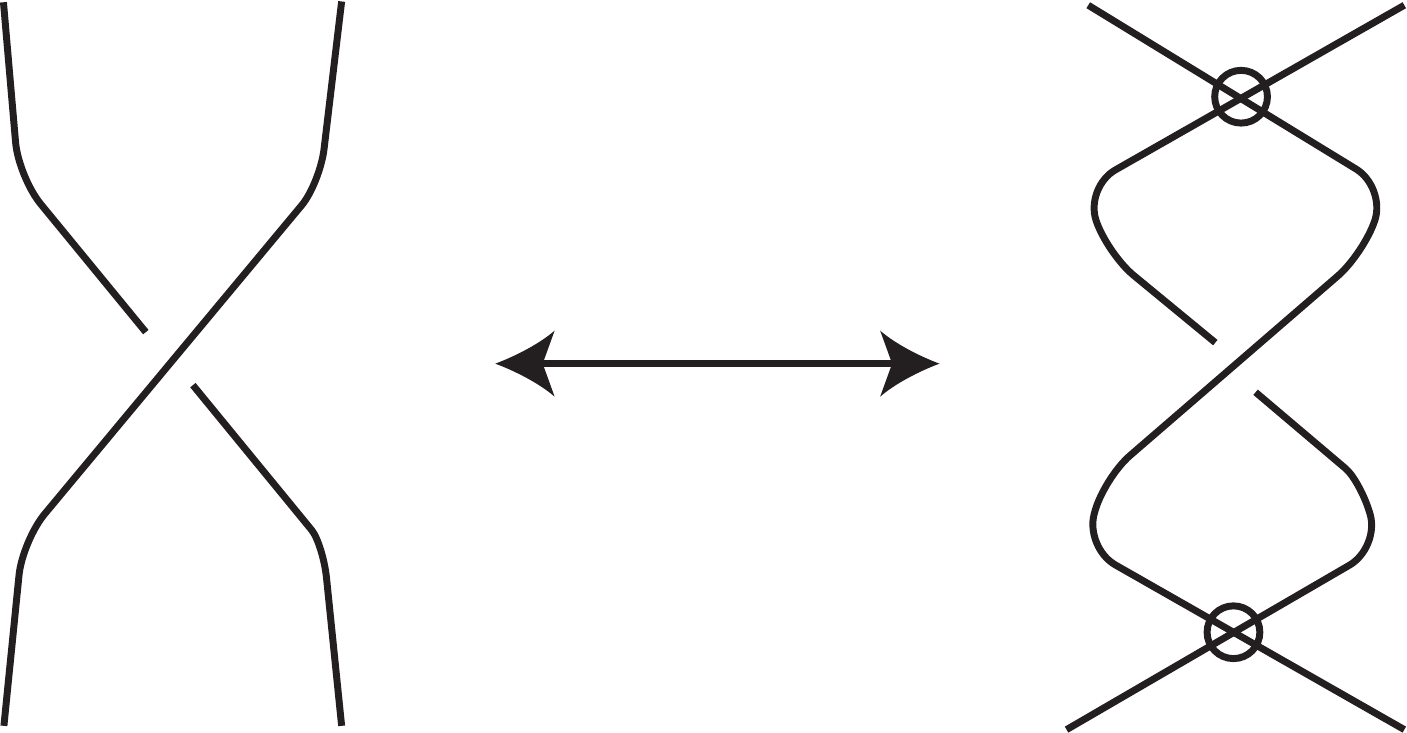}
\caption{Virtualization} \label{twist2}
 \end{figure}

The geometrical sense of Reidemeister moves for free links is that a framed graph is not assumed embedded in any surface. However, when applying a Reidemeister move, one assumes the existence of some ``local'' loop, bigon or triangle. Free knots were first considered by V.\,G.~Turaev \cite{Turaev} who conjectured these knots to be all trivial. 
But later V.~O.~Manturov and then A.~Gibson disproved this conjecture \cite{ManOFK, Ma22, Gib}.
Free knot theory is intimately related to flat virtual knot theory. Let us factorize the theory
of flat virtual knots (or links) by yet another move --- {\em virtualization}, see Fig.~\ref {twist2}, and the new theory was proved to be the free knot theory \cite{Ma22}. 
One may think of a virtualization as a way of changing the immersion
of a framed $4$-graph in the plane such that the cyclic order of
half-edges changes but opposite edges remain opposite.
The exact statement connecting virtual knots and free knots is the following, which easily follows from the definition and has been used by \cite{IMN11}. 
We shall represent free knots by virtual knot diagrams.

\begin{proposition}[\cite{Ma22}]
Two representatives of free links represent the same equivalence class if and only if the corresponding virtual link diagrams are the same modulo a combination of the following transformations:
\begin{enumerate}
\item The generalized Reidemeister moves for virtual knot theory.
\item Crossing switches that make a diagram flat.
\item Virtualization (regardless of the embedding of a representative of a free link).
\end{enumerate}
\end{proposition}

\begin{remark}[\cite{ManOFK}]
The equivalence of free links is coarser than
the equivalence of flat virtual links: our free links do not require
any surface. Every time one applies a Reidemeister move to a framed
$4$-graph, one embeds this graph into a $2$-surface arbitrarily
(with framing preserved), apply this Reidemeister move inside the
surface and then forget the surface again.
\end{remark}

Since free knot theory allows us to have more flexibility on generalized Reidemeister moves, we may use it to study different invariants of virtual knots and
refinements of invariants of knots and cobordisms in higher
dimensions \cite{Cobordisms}. Free links have natural applications to many problems
concerning embedding of graphs into surfaces \cite{Man4, IMN2}.

\section{Towards free braid theory}

In classical knot theory, knots and links can be represented as equivalence classes of braids modulo Markov moves  \cite{LR95, Birman, Markov}.
Let $\{(x,y,z)\in\R^3|y=0,z=0\}$ and $\{(x,y,z)\in\R^3|y=0, z=1\}$ be two parallel lines in $\R^3$ with $n$ points on each of them, with $x$-coordinates $1, 2, \ldots, n$ respectively. 
An {\em $n$-strand braid} is a set of $n$ non-intersecting smooth paths in $\R^3$ connecting the chosen points on the first line with the chosen points on the second line, so that there are no two paths leading to the same point and so that there are no local maxima or minima with respect to the height function, that is, the third coordinate.
Each of these smooth paths is called a {\em strand} of the braid.
Two braids $B_0$ and $B_1$ are equal if they are isotopic to each other, that is, if there exists a continuous family of braids $B_t, t\in [0,1]$ starting at $B_0$ and finishing at $B_1$.
Under the isotopy equivalence the set of all $n$-strand braids forms a group $B_n$, called {\em the Artin $n$-strand braid group}. Here, the operation is to connect the endpoints of the first braid to the corresponding starting points of the second braid and, the unit element of the group is the braid represented by all vertical parallel strands.

\begin{figure}
\centering\includegraphics[width=280pt]{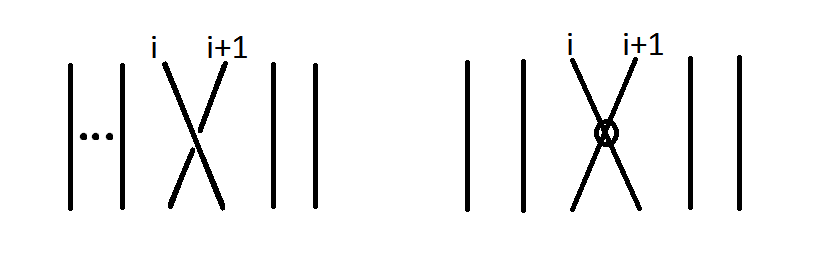} \caption{The classical crossing and the virtual crossing} 
\label{real and virtual crossings}
\end{figure}

It is natural to study braids by using {\em braid diagrams}. 
Given a braid, a diagram is obtained by taking a generic projection of the braid in $\R^{3}$ to a plane parallel to the $xz$-plane.
More precisely, a {\em braid diagram} is a graph lying inside the rectangle $[1, m] \times [0, 1]$ endowed with the following structure and having the following properties:
\begin{itemize} 
\item Points $[0, i] $and $[1, i], i = 1,\ldots, m$ are vertices of valency one, the other points of type $[0, t]$ and $[1, t]$ are not vertices of graph.
\item All other graph vertices (crossings) have valency four; opposite edges at such vertices form angles $\pi$.
\item Unicursal curves, that is, lines consisting of edges of the graph, passing from an edge to the opposite one, go from vertices with first coordinate one and come to vertices with first coordinate zero; they must be descending.
\item Each vertex of valency four is endowed with either an overcrossing or a undercrossing structure.
\end{itemize}
This diagrammatic approach naturally leads to a presentation of the braid group.
The braid diagram can be thought of a composition of the elementary braids $\sigma_i^{\pm} (1\le i\le n-1)$. Here, $\sigma_i$ corresponds the left part of Fig. \ref{real and virtual crossings}.

\begin{definition}
The {\em $n$-strand braid group} $B_n$ is presented by $n-1$ generators $\sigma_1, \ldots, \sigma_{n-1}$ subject to the relations: $\sigma_i\sigma_j=\sigma_j\sigma_i, |i-j|>1$ and $\sigma_i\sigma_{i+1}\sigma_i=\sigma_{i+1}\sigma_i\sigma_{i+1}, 1\le i\le n-2.$
\end{definition}

\begin{remark}
 The latter relations of the presentation of $B_n$ correspond to the third Reidemeister moves applied to braid diagrams on the plane.
The second Reidemeister move appears as the trivial relation $\sigma_{i}\sigma_{i}^{-1}=1$ (or $\sigma_{i}^{-1}\sigma_{i}=1$), and the first Reidemeister move is not applicable to braid diagrams because of the monotonicity of the braid strands.
\end{remark}

\begin{figure}
\centering\includegraphics[width=300pt]{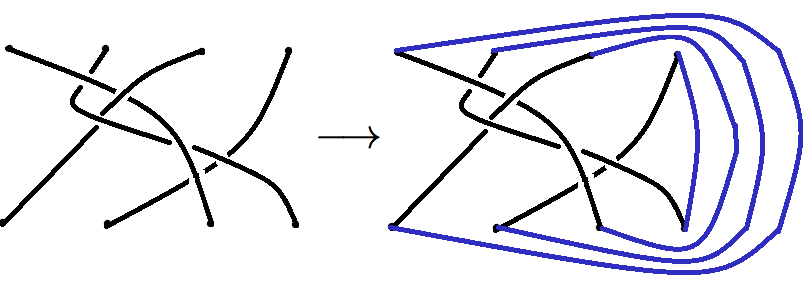} \caption{The closure of a braid} 
\label{closure of a braid}
\end{figure}

With each braid diagram $b$, one can associate a knot (or link) diagram by taking the {\em closure $\Cl(b)$}, obtained by connecting the lower ends of the braid $b$ with the upper ends by simple disjoint arcs, see Fig.\ref{closure of a braid}. 
Obviously, isotopic braids generate isotopic links. 
The inverse process of taking closure is called the {\em braiding} of a knot or a link. The celebrated {\em Alexander theorem} states that for each link $L$ there exists a braid $b$ such that $\Cl(b)$ is isotopic to $L$.
The precise equivalence relation on braids capturing the isotopy of two links as closures of braids is described by the {\em Markov theorem}.
The Markov theorem states that the closures of two braids represent isotopic links if and only if the two braids are related by successive applications of two types of moves, conjugation in the braid groups and Markov stabilization moves, on the set of braids \cite{Birman, Markov}. The Markov theorem is powerful in constructing invariants for classical knots and links \cite{Oht}.

To study invariants for free links, an idea is to investigate the representation of the corresponding braid groups. Hence, it is natural to ask for a Markov type theorem for free links. 
The tool we shall use is virtual braid theory, and a generalized Markov theorem in this setting.
Just as classical braids, virtual braids have a purely algebraic definition using virtual braid diagrams.

\begin{definition}[\cite{Ver}]
An {\em $n$-strand virtual braid diagram} is a graph lying in $\R^2$ with $2n$ vertices of valency one  having coordinates $(i, 0)$ and $(i, 1)$ for $i = 1,\ldots, n$ and a finite number of vertices of degree four. 
The graph is a union of $n$ smooth curves connecting points on the line $\{y = 1\}$ with those on the line $\{y = 0\}$ and descending with respect to the vertical coordinate. The $4$-valent vertices come from intersections. 
Each crossing is either endowed with a structure of overcrossing or undercrossing, as in the case of classical braids, or a virtual crossing, marked by encircling it. See Fig. \ref{virtual braid diagram}.

\begin{figure}
\centering\includegraphics[width=100pt]{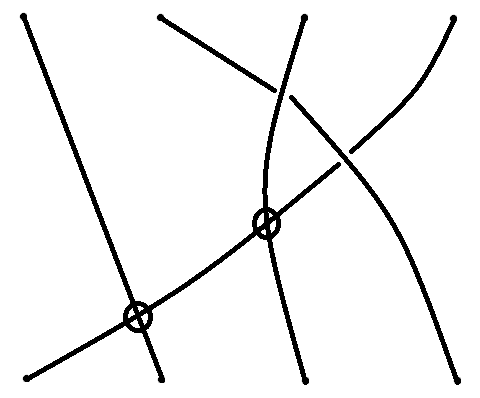} \caption{A virtual braid diagram} 
\label{virtual braid diagram}
\end{figure}

A {\em virtual braid} is an equivalence class of virtual braid diagrams by planar isotopies and all {\em Reidemeister moves for virtual braids}, that is,  the generalized  Reidemeister moves where both the diagram before the move and the diagram after the move are braided. Note that the first classical Reidemeister move and the first virtual Reidemeister move are not in the list. 
\end{definition}

\begin{figure}
\centering\includegraphics[width=300pt]{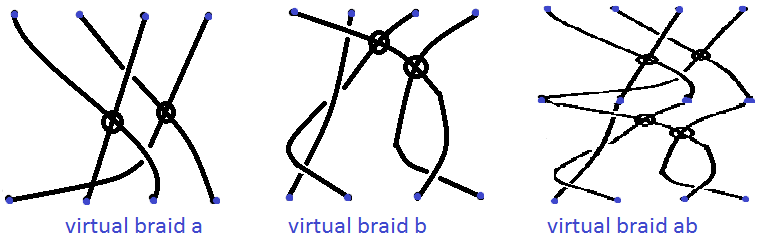} \caption{The product of two virtual braids} 
\label{braidmult}
\end{figure}

\begin{figure}
\centering\includegraphics[width=300pt]{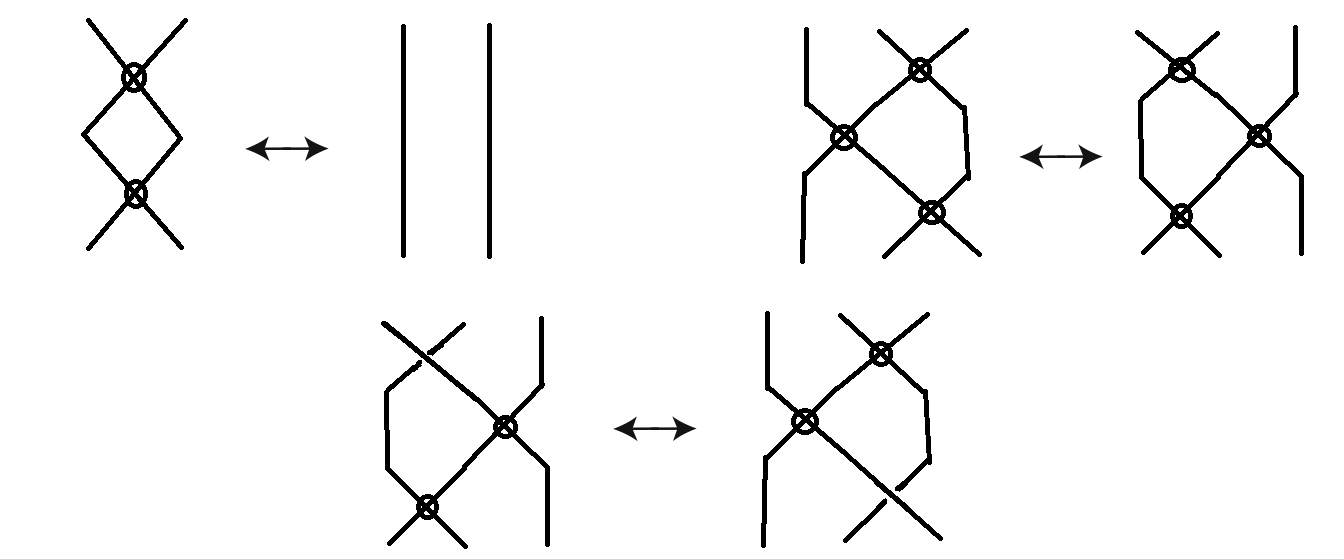} \caption{The braid isotopies} 
\label{braid isotopy 1}
\end{figure}

Like classical braids, $n$-strand virtual braids form a group $vB_n$ with respect to connecting the corresponding end points, smoothing the angles caused by juxtaposition and rescaling the vertical coordinate. See Fig. \ref{braidmult}.
The generators of this group are $\sigma_1,\ldots ,\sigma_{n-1}$ (for classical crossings) and $\zeta_1, \ldots, \zeta_{n-1}$ (for virtual crossings). See Fig. \ref{real and virtual crossings}.
Obviously, beside the classical relations, we have $\zeta_i\zeta_j=\zeta_j\zeta_i$ and $\zeta_i\sigma_j=\sigma_j\zeta_i$ for $|i-j|>1$ by the braid isotopy. Furthermore, we have 
$$\zeta_i^2=1, 1\le i\le n-1,$$ $$\zeta_{i+1}\zeta_i\zeta_{i+1}=\zeta_i\zeta_{i+1}\zeta_i, 1\le i\le n-2,$$ and $$\sigma_i\zeta_{i+1}\zeta_i=\zeta_{i+1}\zeta_i\sigma_{i+1}, 1\le i\le n-2,$$ 
following from the second, the third virtual Reidemeister move and the semi-virtual move respectively. See Fig. \ref{braid isotopy 1}.
The group $vB_n$ given by such a presentation coincides with the group of $n$-stand virtual braids \cite{Ver}.

\begin{remark}
Without loss of generality, we shall assume our braids as polygonal ones, that is, each strand does not have to be smooth, but angles may appear in it. 
\end{remark}

\begin{remark}
A {\em flat virtual braid diagram} is obtained from a virtual braid diagram by replacing all the classical crossings with flat crossings.
We can similarly define the group of the {\em $n$-strand flat virtual braids} by adding the relations $\sigma_i^2=1, 1\le i\le n-1$ to $vB_n$ \cite{KL}. 
\end{remark}

\begin{figure}
\centering\includegraphics[width=300pt]{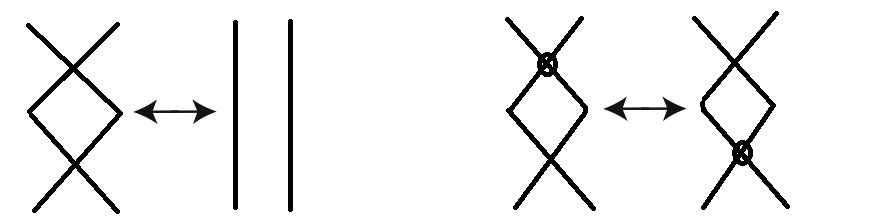} \caption{The free braid isotopies} 
\label{braid isotopy 2}
\end{figure}

Analogously to the definition of free knots and links, we define the {\em $n$-strand free braid group} $fB_n$ as the quotient of the $n$-strand virtual braid group $vB_n$ by two more relations: the cross-switching and virtualization. 
In principal, we add two extra relations to the presentation of the group $vB_n$: for any $1\le i\le n-1$, 
$$\sigma_i^2=1\qquad\text{ and }\qquad\sigma_i\zeta_i=\zeta_i\sigma_i.$$ See Fig. \ref{braid isotopy 2}. Equivalently, $fB_n$ is an equivalence class of $n$-strand flat braid diagrams by isotopies, virtualization and Reidemeister moves for flat virtual braid diagrams (not including the first Reidemeister move and the first virtual Reidemeister move). 
The algebraic definition of the free braid group is summarized as follows.

\begin{definition}
The set of {\em the $n$-strand free braids} $\mathrm{fB}_n$ is a group with $2n-2$ generators $\sigma_1, \sigma_2, \ldots, \sigma_{n-1}, \zeta_1, \zeta_2, \ldots, \zeta_{n-1}$ subject to the following relations: 
\begin{itemize}
\item (Relations for classical braids) 
\begin{itemize}
\item $\sigma_i\sigma_j=\sigma_j\sigma_i$, for all $|i-j|>1, 1\le i, j\le n-1$; 
\item $\sigma_{i}\sigma_{i+1}\sigma_i=\sigma_{i+1}\sigma_i\sigma_{i+1}$, for $1\le i\le n-2$.
\end{itemize}
\item (Additional relations for virtual braids) 
\begin{itemize}
\item $\zeta_i\zeta_j=\zeta_j\zeta_i$ and $\zeta_i\sigma_j=\sigma_j\zeta_i$, for all $|i-j|>1$; 
\item $\zeta_i\zeta_{i+1}\zeta_i=\zeta_{i+1}\zeta_i\zeta_{i+1}$ and $\sigma_i\zeta_{i+1}\zeta_i=\zeta_{i+1}\zeta_i\sigma_{i+1}$ for $1\le i\le n-2$; 
\item $\zeta_i^2=1$.
\end{itemize} 
\item (Additional relations for free braids) $\sigma_i\zeta_i=\zeta_i\sigma_i$, $\sigma_i^2=1$ for all $1\le i\le n-1.$
\end{itemize}
\end{definition}

\begin{remark}
To the best of our knowledge, the notion of free braid has not appeared in the literature; possibly, this was because the existence of non-trivial  free knots was proved not so much time ago.
\end{remark}

In \cite{KL2} Louis Kauffman and Sofia Lambropoulou proved  the Alexander's theorem for virtual links, that is, every virtual link is isotopic to the closure of some virtual braid. The Alexander for free links follows easily as a corollary. 

\begin{theorem}\label{AlexanderFree}
For any free link, there exists a free braid whose closure is isotopic to the given link.
\end{theorem}

\begin{proof}

By the natural surjective map from virtual knot theory to free knot theory $\phi: \{\text{Virtual links}\}\rightarrow\{\text{Free links}\}$, a free link can be represented by a flat virtual link, so we apply to it a known braiding algorithm from  \cite{KL, KL2}, and so the theorem follows, assuming the Alexander theorem for flat virtual links. 
\end{proof}

To show Markov theorem for free links, the proof of Theorem \ref{AlexanderFree} is not sufficient. We need the generalization of the braiding algorithm presented in \cite{KL2} to the case of free links.
According to \cite{KL2}, the following assumptions are made on oriented virtual links:

\begin{itemize}
\item Every virtual link diagram is piecewise linear, that is, a union of line segments called {\em arcs}. For such polygonal links \cite{LR95}, there is another type of ``moves": {\em subdivision} of an arc into two smaller arcs by marking it with a point \cite{KL2}.
\item By locally isotopic shifts, a virtual link is supposed to be in {\em general position}, where the diagram consists of intervals and satisfies following conditions: 
\begin{itemize}
\item There are no horizontal arcs in the diagram.
\item On the horizontal and vertical level of a small neighborhood of each crossing, there are no other crossings or subdividing points.
\item When zooming into a small neighborhood of a crossing, the two involved arcs are either going up or going down.
\end{itemize}
\end{itemize}

We shall present an algorithm for constructing a polygonal free braid from a free link.
The idea of the algorithm in \cite{KL2} is summarized as follows. 
All arcs in the diagram are oriented either upwards or downwards. We want to obtain a braid diagram where all the arcs are oriented downwards. This process is called {\em braiding}. The braiding algorithm for a virtual braid diagram in general position is to eliminate all arcs oriented upwards, called {\em up-arcs}. 
For each up-arc we  cut at a point near its upper end and pull the upper half upward and the lower half downward, and creating a pair of new braid strands, by using isotopies and braid moves. Here the endpoints of the new strands are exactly the cut points. All the new crossings caused by adding these new strands are assumed to be virtual, which ensures that the closure of this new braid is the same as the closure of the initial braid. 
It does not matter in which order the elimination of up-arcs happen. 

\begin{remark}
Since free links are virtual links modulo two equivalence relations, we may use piecewise flat virtual links to represent the free links, and apply the braiding algorithm for flat virtual links on polygonal free links in general position. The resulting free braid is simpler as we have more equivalence relations.
For example, the braiding of a local crossing consisting of two up-arcs (appeared in \cite{KL2}) is simpler because of the virtualization for free braids. See Fig. \ref{braid2uparcs}.
\begin{figure}
\centering\includegraphics[width=250pt]{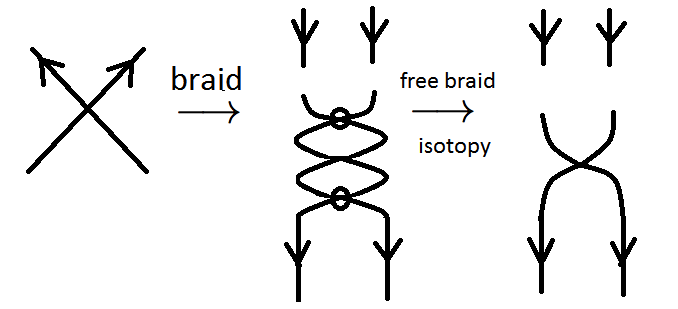} \caption{Braiding the crossing of two up-arcs in a free link} 
\label{braid2uparcs}
\end{figure}
An example of braiding for a free link, using the algorithm given in \cite{KL2}, is shown in Fig. \ref{braiding}. 
\begin{figure}
\centering\includegraphics[width=350pt]{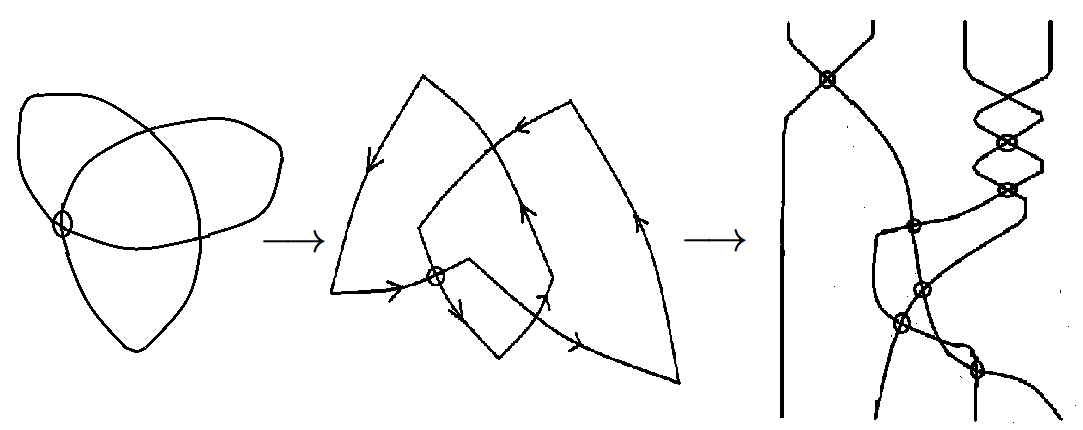} \caption{The process of braiding for a free knot} 
\label{braiding}
\end{figure}
\end{remark}

\section{$L$-moves and Markov theorem}

Let us first recall the generalized Markov theorem for virtual knot theory, and then, apply it to free links. Virtual braids close up into virtual link diagrams. Obviously, isotopic virtual braids close to isotopic virtual links. Furthermore, all virtual link isotopy classes can be represented by closures of virtual braids. In \cite{Kam}, Seiichi Kamada proved an analogue of Markov theorem for the case of virtual braids. In \cite{KL2} Louis Kauffman and Sofia Lambropoulou using the $L$-move method to give a local version of the Markov theorem for virtual braids. This method provides a one-move Markov theorem in classical braid theory and applies to many diagrammatic settings related to virtual braids. 
We emphasize that the Markov theorem for flat virtual links follows immediately from the argument of \cite{KL2}, but for free braids it is not so: we need to deal with the virtualization move. We formulate the definition of $L$-move for a free braid as follows.

\begin{definition}[\cite{KL2}]\label{Lmove}
An {\em $L$-move} on a free braid is a move on the flat virtual braid diagram, representing the free braid, obtained by cutting open an arc of the diagram, and pulling the two ends of the gap, so as to create a new pair of braid strands, whose closure is isotopic to the closure of the original flat virtual braid and which intersects virtually other parts of the braid diagram.

\begin{enumerate}
\item If the new pair of strands is at the cutpoint, the $L$-move is called the {\em basic $L$-move}. 
See Fig. \ref{L-move} from (a) to (b).

\begin{figure}
\centering\includegraphics[width=350pt]{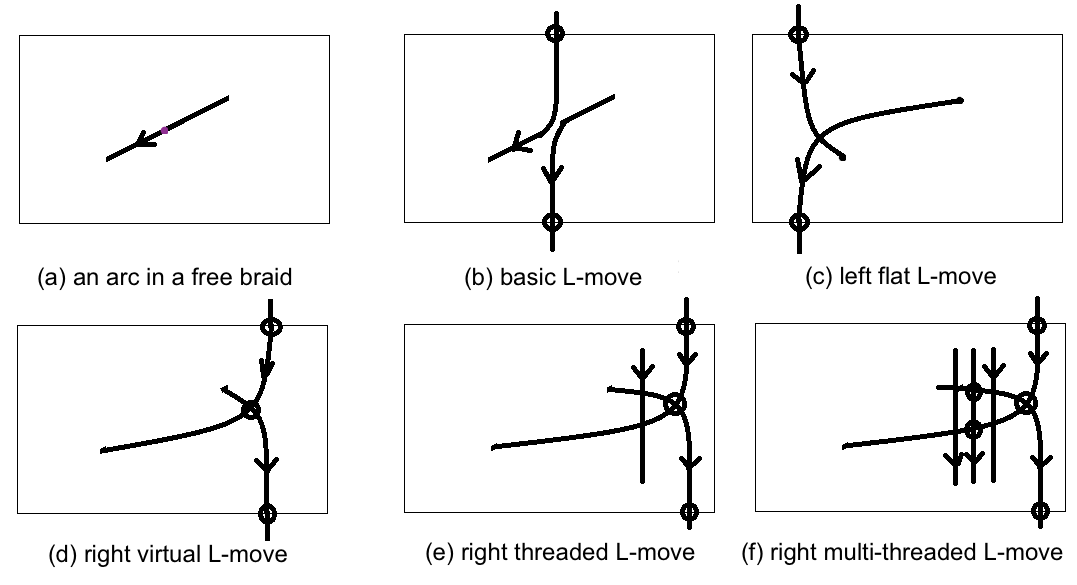} \caption{Examples of $L$-moves of an arc in a flat virtual braid diagram} 
\label{L-move}
\end{figure}

\item If the new pair of strands is to the right of the cutpoint and if there are no other stands between the cutting point and the new strands, there will be a new crossing created by the new strands. Depending on the type of the crossing, flat or virtual, the $L$-move is called {\em right flat $L$-move or right virtual $L$-move}. And if we change the "right" to the "left" in the above sentence, the move is called {\em left flat $L$-move or left virtual $L$-move}.
See Fig. \ref{L-move} from (a) to (c) or (d).

\item If the new strands are to the left/right of the cutpoint, with a new virtual crossing coming from the intersection of the new strands,  and if there is one more strand in between the cutpoint and the new pair of strands, with flat crossings, the move is called {\em left/right threaded $L$-move}.
See Fig. \ref{L-move} from (a) to (e).
\end{enumerate}
\end{definition}

\begin{remark}
In \cite{KL2}, more $L$-moves were used in the proof of the Markov theorem for virtual links. For example, {\em multi-threaded $L$-moves}, where there are more than one strand between the cutting point and the new strands (Fig. \ref{L-move} (f)), is used. Such $L$-moves can be reduced to a composition of other $L$-moves and free braid isotopy according to \cite{KL2}. So they do not appear in the statement of the theorem. The same happens with the statement of the  Markov theorem for free links. 
\end{remark}

The Markov theorem for free links is formulated as follows.

\begin{theorem}\label{freemarkov}
Two oriented free links are isotopic if and only if two corresponding free braids differ by a finite sequence of free braid isotopy and the following moves and their inverses:
\begin{enumerate}
\item Flat conjugation.
\item Right virtual $L$-moves.
\item Right flat $L$-moves.
\item Right and left threaded $L$-moves. 
\end{enumerate}
\end{theorem}

\begin{figure}
\centering\includegraphics[width=250pt]{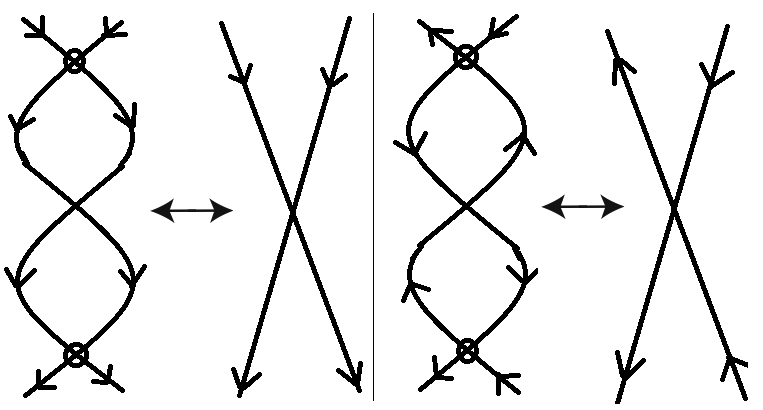} \caption{Local picture of two free links which differ by a virtualization} 
\label{virtualizatininlink}
\end{figure}

First of all, we essentially need to work on one virtualization move, shown in the left part of Fig. \ref{virtualizatininlink}. More precisely, the following lemma holds.

\begin{figure}
\centering\includegraphics[width=250pt]{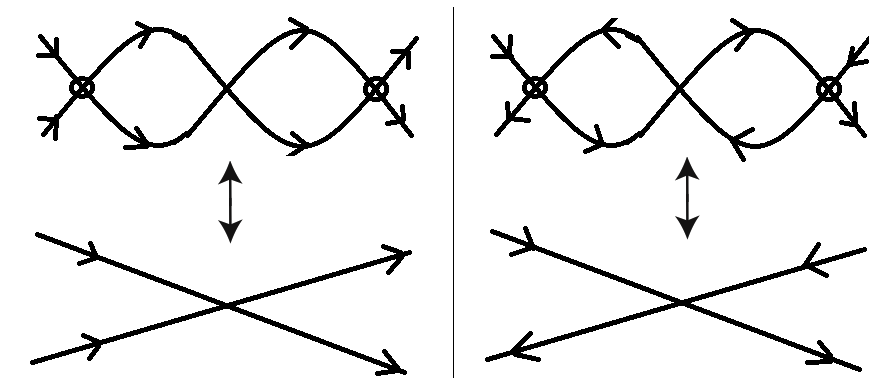} \caption{Rotating Fig. \ref{virtualizatininlink} counterclockwise by 90 degrees} 
\label{virtualizatininlinks}
\end{figure}

\begin{lemma}\label{key}
\begin{enumerate}
\item The left illustration of Fig. \ref{virtualizatininlink} and the upper right illustration of Fig. \ref{virtualizatininlinks} differ by detour moves.
\item Respectively, the horizontal mirror of the left illustration of Fig. \ref{virtualizatininlink} and the horizontal mirror of the upper right illustration of Fig. \ref{virtualizatininlinks} differ by detour moves.
\item The right illustration of Fig. \ref{virtualizatininlink} and the vertical mirror of the upper left illustration of Fig. \ref{virtualizatininlinks} differ by detour moves.
\item Respectively, the vertical mirror of the right illustration of Fig. \ref{virtualizatininlink} and the upper left illustration of Fig. \ref{virtualizatininlinks} differ by detour moves.
\end{enumerate}
\end{lemma}
\begin{proof}
Two corresponding pieces of virtual link diagrams differ by two detour moves. See Fig.~\ref{proof5}. The proof of the first statement is accomplished if we replace all classical crossings in the figure by flat crossings. The rests follow from appropriately changing the arrows in the diagram.
\begin{figure}
\centering\includegraphics[width=330pt]{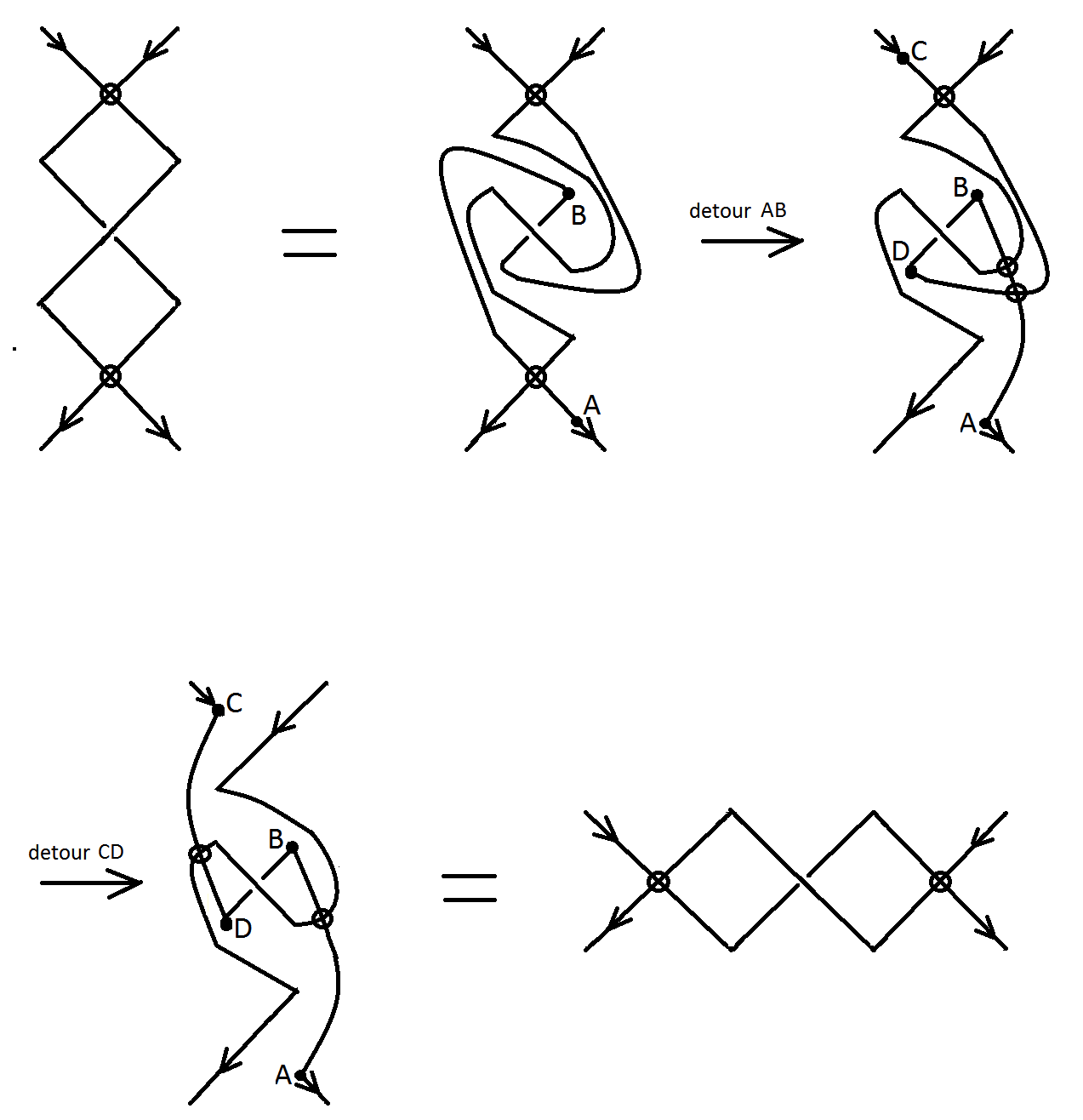} \caption{Proof of Lemma \ref{key}} 
\label{proof5}
\end{figure} 
\end{proof}

\begin{remark}
The only thing we need to do in order to prove Theorem \ref{freemarkov} which does not follow immediately from \cite{KL2}, is that the virtualization move can enter the game: whenever two free link diagrams differ by a virtualization, we may arrange it so that the corresponding pieces are ``well braided" so that there is no need to perform any further effort.
Lemma \ref{key} plays a key role in the proof of our main theorem. Namely, whenever we apply a virtualization move, we may always assume that it is done in the preferred way, which locally agrees with the braiding. 
\end{remark}

\begin{proof}[Proof of Theorem \ref{freemarkov}]

By Theorem 4 \cite{KL2}, a corollary of the main theorem in the paper, two flat virtual links are isotopic if and only if their corresponding flat virtual braids are connected by a finite sequence of flat braid isotopies and the moves listed in the theorem.
Two free links $L_1$ and $L_2$ are isotopic if and only if $L_1$ can be transformed into $L_2$ by the isotopy on flat virtual links and virtualization. 
So we only need to prove that if two free links differ by a virtualization, then their corresponding free braids are connected by the moves listed in the statement of the theorem. 
Assume we have two oriented free link diagrams in general position, which are identical except for a virtualization, with all possible orientations of arcs.
In the braiding process, all the new crossings as a result of the new added braid strands are going to be virtual. 
Hence, using the free braid isotopy, we only need to compare the braiding of the corresponding different local pieces in two free link diagrams. We place the two free link diagrams in $\R^2$ so that:
\begin{enumerate}
\item For a fixed orientation, the two corresponding pieces fall into either of the two cases in Fig. \ref{virtualizatininlink}. 
\item In the opposite orientation of the strands, the two corresponding pieces fall into either the horizontal mirror image of the left two illustrations in Fig. \ref{virtualizatininlink} or the vertical mirror of the two illustrations in Fig. \ref{virtualizatininlink}.
\end{enumerate}
In view of Lemma \ref{key}, it is sufficient to consider only (1).

In the first case of (1), which corresponds to the left two pictures of Fig. \ref{virtualizatininlink}, the two pieces in the corresponding free link diagram stay unchanged during the braiding process. Therefore by the defining relation $\zeta_i\sigma_i\zeta_i=\sigma_i$ in the free braid, the two free braids corresponding to the two free links represent the same element in the free braid group.

In the second case of (1) which corresponds to the picture on the right hand side of Fig. \ref{virtualizatininlink}, by Lemma \ref{key}, this case is reduced to the first case because the braids corresponds to detour moves and rotation of links are connected by the $L$-moves and the flat conjugation (As a consequence of Markov theorem for flat virtual links \cite{KL2}, the braiding of an oriented flat virtual link and the braiding of its rotation by $90$ degrees in $\R^2$ are connected by $L$-moves and conjugation). 
Now, the proof of Theorem \ref{freemarkov} is concluded.
\end{proof}


\begin{remark}
There is also a type of knot theory similar to free links, defined by virtual knots modulo virtualization. We call them {\em $V/Z$ links}. Combining the proof of Lemma \ref{key} and the similar argument of  the proof of Theorem \ref{freemarkov} and the arguments in \cite{KL2}, we obtain the following Markov theorem for the $V/Z$ links. The only difference in the proof is that we need to take care of the over crossings and the undercrossings. To state the theorem, one has to use the terminology in \cite{KL2}, that is, $L$-moves for virtual links.
\end{remark}

\begin{theorem}
Two oriented $V/Z$ links are isotopic if and only if two corresponding $V/Z$ braids differ by virtual braid isotopy, virtualization and a finite sequence of the following moves and their inverses:
\begin{enumerate}
\item Real conjugation.
\item Right virtual $L$-moves.
\item Right real $L$-moves.
\item Right and left under threaded $L$-moves. 
\end{enumerate}
\end{theorem}

\section{Further remarks}

The main motivation for our interest to Markov's theorem for free links is that we would like to construct invariants for free knots and links out of free braids.  
Up to now, all previously known invariants of free knots and links are in some sense "perpendicular" to classical or virtual invariants; they were based on the notion of parity \cite{IMN11} and not on the gadgetry known in the classical case. We are going to start tackling problems concerning free link invariants by using classical objects, such as the Yang-Baxter equation.
We have a presentation of a free braid group and we would like to know how it helps in the classification of the free links. We want to start with investigating simple questions such as, what free braid diagrams have trivial knot/link closures? For example, the closure of braid $\sigma_1\xi_1$ is not trivial.

In classical knot theory, there are many standard ways of constructing invariants of knots and links. For example, the {\em quantum invariants} (cf. for example \cite{Oht} for the precise definition and the connection to statistical mechanics). 
Let $V$ be a vector space of dimension $n$ and $R: V\otimes V\rightarrow V\otimes V$ be a linear transformation, that is, $R$ is a matrix of size $n^2$. We denote by $R^{a b}_{c d}$ the $(a\otimes b, c\otimes d)$-th matrix entry, where $a, b, c$ and $d$ belong to a set of the basis of $V$. 

\begin{definition}
Let $\mathrm{Id}$ be the identity map from $V$ to itself, then {\em the Yang-Baxter equation} is an equation on $V^{\otimes3}$ given by
\begin{equation}\label{YBE}
(R\otimes \mathrm{Id})(\mathrm{Id}\otimes R)(R\otimes \mathrm{Id})=(\mathrm{Id}\otimes R)(R\otimes \mathrm{Id})(\mathrm{Id}\otimes R).
\end{equation}
\end{definition}

\begin{remark}
Note that the Yang-Baxter relation here is considered over some field $\mathbb{K}$, that is, $V$ is a vector space over a field $\mathbb{K}$. 
However, a similar calculation of solutions to the Yang-Baxter equation can be performed on a module $V$, over a non-commutative ring $\mathbb{K}$ or a ring $\mathbb{K}$ with zero divisors.
\end{remark}

Let us take the following representation of the braid group $B_n$ on $V^{\otimes n}$
\begin{equation}
\sigma_{i}\to \mathrm{Id}\otimes\dots \otimes\mathrm{Id}\otimes \underbrace{R}_{i\text{th,} (i+1) \text{st}}\otimes\mathrm{Id}\otimes\dots\otimes \mathrm{Id},
\end{equation}
where $R$ corresponds to the $i$-th and $(i+1)$-st factors in $V^{\otimes n}$.
Thus we obtain a representation of $B_n$. 
In particular, the relation $\sigma_i\sigma_{i+1}\sigma_{i}=\sigma_{i+1}\sigma_{1}\sigma_{i+1}$ corresponds to the Yang-Baxter equation (\ref{YBE}).
The first step to obtain an invariant for knots and links, one can take
the trace of this representation, and this trace is automatically invariant under conjugation (which is known as the first Markov move \cite{Ka2, Man1}).
Nevertheless, the problem of finding such representations for which the trace is invariant under
all Markov moves is rather complicated, see, for example, \cite{Oht}.

For virtual knot theory this problem becomes even more complicated because we have two sorts
of ``stabilization moves" corresponding to the first classical Reidemeister move and the first virtual
Reidemeister move.
So, in the present paper we restrict ourselves to a couple of examples.

If we want to study the invariants of free knots by investigating solutions for Yang-Baxter equation (\ref{YBE}). Notice that there are two extra relations imposed on $R$: 
\begin{equation}\label{Rvirtualization}
  R^{ab}_{cd}=R^{ba}_{dc}\text{, for all }a, b, c, d \text{ being in the set of the basis of } V;  
\end{equation}
\begin{equation}\label{Rflat}
R^{2}=I_{n_2} \text{ ($I_{n_2}$ is the identity matrix of size $n^2$)}. 
\end{equation}
Here, (\ref{Rvirtualization}) comes from the version of the virtualization $\sigma_i\xi_i=\xi_i\sigma_i$ and (\ref{Rflat}) comes form the relations $\sigma_i^2=1$ and $\xi^2=1.$
Each solution of the set of equations (\ref{YBE}), (\ref{Rvirtualization}) and (\ref{Rflat}) gives rise to a representations of the free braid group $fB_n$.

\begin{example}
Let $V$ be a $2$-dimensional vector space/module over a field/ring $\mathbb{K}$, with a basis $\{e_0, e_1\}$. A linear map $R\in\mathrm{End}(V\otimes V)$ is presented by $R=\begin{pmatrix} R_{00}^{00} & R_{00}^{01} & R_{00}^{10} & R_{00}^{11} \\ R_{01}^{00} & R_{01}^{01} & R_{01}^{10} & R_{01}^{11} \\ R_{10}^{00} & R_{10}^{01} & R_{10}^{10} & R_{10}^{11} \\ R_{11}^{00} & R_{11}^{01} & R_{11}^{10} & R_{11}^{11} \end{pmatrix}$ and we want to solve $R$ satisfying the set of equations (\ref{YBE}), (\ref{Rvirtualization}) and (\ref{Rflat}). 
In order not to deal with the general case with too many equations, we shall restrict ourselves to the case with many zero elements, the so-called eight-vertex model, that is, $R=\begin{pmatrix} a & 0 & 0 & b \\ 0 & c & d & 0 \\ 0 & d & c & 0 \\ b & 0 & 0 & a \end{pmatrix}$. Here $R$ satisfies the equation (\ref{Rvirtualization}).
Then by solving equations (\ref{YBE}) and (\ref{Rflat}) we obtain the following set of equations 
\begin{align*}
 a^2b+bca&=abc+bd^2; \\
 cd^2+dac&=cda+db^2; \\
 dc^2+cad&=adc+b^2d; \\
 dcd+cac&=aca+bab; \\ 
 ad^2+bc^2&=d^2a+cb^2; \\
 d^2b+cba&=acb+ba^2; \\
 a^2+b^2&=c^2+d^2=1; \\
 ab+ba&=cd+dc=0.
 \end{align*}
If furthermore, $\mathbb{K}$ is commutative with no zero divisors, then the set of equations is reduced to 
\begin{align*}
b^2d=bd^2&=2ab=2cd=0; \\
ac^2&=a^2c; \\
bc^2&=b^2c; \\
a^2+b^2&=c^2+d^2=1.
\end{align*}
Hence, we obtain the complete set of solutions of $(a, b, c, d)$ by solving the set of equations: $(1, 0, 0, \pm1),$ $(1, 0, 1, 0),$ $(-1, 0, 0, \pm1),$ $(-1, 0, -1, 0),$ $(0, 1, 1, 0)$ and $(0, -1, -1, 0).$

If $\mathbb{K}$ is a ring with zero divisors, then there are more interesting solutions. For example, when $\mathbb{K}=\Z_{12}$, a solution to $(a,b,c,d)$ is $(4, 3, 3, 4).$
\end{example}

\begin{remark}
It is interesting to study the representations of free braid groups when $\mathbb{K}$ as a ring has zero divisors or is noncommutative. 
Furthermore, it is much more complicated to study quantum invariants for free knots and links, and we will investigate this problem in a subsequent paper.
\end{remark}

\begin{remark}
Other ways of constructing virtual knot invariants are due to A. Bartholomew, R. Fenn and L. Kauffman, V.O.Manturov and they consist of constructing biquandles \cite{BF06, KM}. In this biquandle setup instead of associating vector spaces to strands and their tensor powers to link diagrams, we associate an $n$-dimensional space to an $n$-strand braid and write a ``simplified analogue of the Yang-Baxter equation". We shall also touch on all such questions in a subsequent paper. 
\end{remark}

\begin{remark}
We still do not know whether free knots and links are algorithmically recognizable, unlike virtual knots and links.
However, free braids are a much simpler object since they form a group. Possibly, there is an algebraic way to recognize free braids using the method of Bardakov \cite{Bar}.
\end{remark}

{\bf Aknowledgement:} We would like to express our genuine gratitude to Louis Kauffman and Sofia Lambropoulou for their fruitful discussions and helpful comments.

\end{document}